\documentclass[12pt]{article}
\usepackage{amsmath}
\usepackage{amssymb}
\usepackage{amscd}
\usepackage{amsthm}

\theoremstyle{plain}
\newtheorem{Thm}{Theorem}[section]
\newtheorem{Conj}[Thm]{Conjecture}
\newtheorem{Prop}[Thm]{Proposition}
\newtheorem{Cor}[Thm]{Corollary}
\newtheorem{Lem}[Thm]{Lemma}

\theoremstyle{definition}

\newtheorem{Rem}[Thm]{Remark}

\numberwithin{equation}{section}

\title{Derived equivalence for stratified Mukai flop on $G(2,4)$}
\author{Yujiro Kawamata}

\begin{document}

\maketitle

\section{Introduction}

We consider a conjecture on the derived equivalence for 
$K$-equivalent varieties in the case of a certain flop between $9$-dimensional 
varieties. 
The following is a very general version of this conjecture 
on the fully faithful derived embedding for log $K$-related varieties 
(\cite{Kawamata3}~Conjecture~2.2):

\begin{Conj}
Let $(X, B)$ and $(Y, C)$ be pairs of quasi-projective varieties 
with $\mathbb{Q}$-divisors such that
there exist quasi-finite and surjective morphisms 
$\pi: U \to X$ and $\sigma: V \to Y$ from smooth varieties, 
which may be reducible, 
such that $\pi^*(K_X + B) = K_U$ and $\sigma^*(K_Y+C) = K_V$.
Let $\mathcal{X} \to X$ and $\mathcal{Y} \to Y$ be natural morphisms from
the associated Deligne-Mumford stacks. 
Assume that there are proper
birational morphisms $\mu: W \to X$ and $\nu: W \to Y$ 
from a third variety $W$ such that $\mu^*(K_X+B) \le \nu^*(K_Y+C)$.
Then there exists a fully faithful exact functor 
$D^b(\text{Coh}(\mathcal{X})) \to D^b(\text{Coh}(\mathcal{Y}))$.
\end{Conj}

The conjecture is proved to be true in some cases in 
\cite{BO}, \cite{BKR}, \cite{Bridgeland2}, 
\cite{Kawamata1}, \cite{Chen}, 
\cite{Namikawa1}, \cite{Kawamata2}, \cite{Kawamata3} and \cite{BK}.

On the other hand, Namikawa \cite{Namikawa2}
proved that a naturally defined functor between the derived categories for 
the stratified Mukai flop on $G(2,4)$ is not 
an equivalence.
This is a flop between $9$-dimensional varieties which is the total space of a 
$3$-parameter degeneration of standard $3$-dimensional flops of 
$(-1,-1)$-curves and the fiber over the most degenerate point is isomorphic 
to $G(2,4)$.

So it is worthwhile to check the conjecture in this special case.
We shall prove that there is nevertheless another functor between the 
same categories which is an equivalence.

%%%%%%%%%%%%%%%%%%%%%%%%%%%%%%%%%%%%%%%%%%%%%%%%%%%%%%%%%%%%%%%%%%%%

\section{Stratified Mukai flop}

We recall the construction of stratified Mukai flops due to Markman 
\cite{Markman} (see also \cite{Namikawa2}).

Let $G = G(r,n)$ be the Grassmann variety of $r$-dimensional subspaces in an
$n$-dimensional vector space $V$.
We assume that $2r \le n$.
Let $S$ (resp. $Q$) be the tautological subbundle (resp. quotient) 
bundle on $G$:
\[
0 \to S \to V \otimes \mathcal{O}_G \to Q \to 0.
\]
The polarization of $G$ is given by 
\[
\text{det }Q = \mathcal{O}_G(1).
\]
We have $\Omega^1_G \cong \mathcal{H}om(Q,S)$.
Since 
\[
\text{Ext}^1(\mathcal{O}_G, \Omega^1_G) 
\cong \text{Ext}^1(Q, S) \cong \mathbb{C}
\] 
there exists a non-trivial extension of vector bundles on $G$ induced by 
the natural homomorphism $\mathcal{O}_G \to \mathcal{H}om(S,S)$:
\[
\begin{CD}
0 @>>> \Omega^1_G @>>> \tilde{\Omega}^1_G @>>> \mathcal{O}_G @>>> 0 \\
@. @| @VVV @VVV \\
0 @>>> \mathcal{H}om(Q,S) @>>> \mathcal{H}om(V \otimes \mathcal{O}_G,S) @>>> 
\mathcal{H}om(S,S) @>>> 0.
\end{CD}
\]

Let $X_0 = T^*G$ and $X$ be the total spaces of $\Omega^1_G$ and 
$\tilde{\Omega}^1_G$, respectively, and $\pi: X \to G$ the projection.
We write $\pi^*S = S_X$ and $\pi^*\mathcal{O}_G(1) = \mathcal{O}_X(1)$.
A point of $X$ is given by a pair consisting of a 
point $p \in G$ and a homomorphism $A: V \to S_p \subset V$ 
which induces a homothety $t: S_p \to S_p$ ($t \in \mathbb{C}$).
We have $(p, A) \in X_0$ if and only if $t=0$.
Let 
\[
Z = \{A \in \text{End}(V) \mid \exists t \in \mathbb{C} \text{ s.t. } 
A^2 = tA, \text{rank }A \le r\}
\]
with a morphism $t: Z \to \mathbb{C}$, and set $Z_0 = t^{-1}(0)$.
Then there is a natural morphism $\phi: X \to Z$
given by $\phi(p, A) = A$, which induces a morphism $\phi_0: X_0 \to Z_0$.
We know that $\phi_0$ and $\phi$ are projective birational morphisms.
We write 
\[
W^{(k)} = \{(p, A) \in X \mid \text{rank }A \le k\}.
\]
Then the exceptional loci of both $\phi_0$ and $\phi$ are equal to 
$W = W^{(r-1)}$.
We have $\dim G = r(n-r)$, $\dim X_0 = 2r(n-r)$, $\dim X = 2r(n-r)+1$, and
$\dim W = 2(r-1)(n-r+1)+n-2r+1$.

We have a dual construction starting from the dual vector space $V^*$.
Let $G^+$ be the Grassmann variety of $r$-dimensional subspaces in $V^*$,
$S^+$ (resp. $Q^+$) the tautological subbundle (resp. quotient bundle) 
on $G^+$, $\text{det }Q^+ = \mathcal{O}_{G^+}(1)$,
$X_0^+$ (resp. $X^+$) the total space of 
$\Omega^1_{G^+}$ (resp. $\tilde{\Omega}^1_{G^+}$), 
$\pi^+: X^+ \to G^+$ the projection, 
$\pi^{+*}S^+ = S^+_{X^+}$, and $\pi^{+*}\mathcal{O}_{G^+}(1) 
= \mathcal{O}_{X^+}(1)$.
There is a natural projective birational morphism $\phi^+: X^+ \to Z$
given by $\phi^+(q, B) = {}^tB$ which induces a 
projective birational morphism $\phi_0^+: X_0^+ \to Z_0$.
We write 
\[
W^{+(k)} = \{(q, B) \in X^+ \mid \text{rank }B \le k\}.
\]
The exceptional locus of $\phi$ is $W^+=W^{+(r-1)}$.
The diagram 
\[
\begin{CD}
X @>{\phi}>> Z @<{\phi^+}<< X^+
\end{CD}
\]
thus obtained is called a {\em stratified Mukai flop}.
By restricting to the subspaces defined by $t = 0$, 
we have a smaller diagram 
\[
\begin{CD}
X_0 @>{\phi}>> Z_0 @<{\phi^+}<< X^+_0
\end{CD}
\]
which is also called a stratified Mukai flop.

If $r = 1$, then the above diagrams are reduced to a {\em standard flop} and a
{\em Mukai flop}.

The birational map $(\phi^+)^{-1} \circ \phi: X -\to X^+$ is decomposed into 
blow-ups and downs (\cite{Markman}).
We consider only the case $r=2$ in the following.
Let $f_1: X_1 \to X$ be the blowing up along the center $G = W^{(0)}$,
the $0$-section of the projection $\pi$.
Then the strict transform $W'$ of $W = W^{(1)}$ is smooth.
Indeed, the projection $W' \to G$ is smooth and its fibers are 
isomorphic to the cone over 
$\mathbb{P}^1 \times \mathbb{P}^{n-3} \subset \mathbb{P}^{2n-5}$.
Let $f_2: X_2 \to X_1$ be the blowing up along the center $W'$.
Let $E_i$ be the exceptional divisor of $f_i$ for $i=1,2$, and 
$E_1'=f_{2*}^{-1}E_1$, the strict transform of $E_1$.
We set $f = f_2 \circ f_1$.
If $r \ge 3$, then we have similar construction with more blow-ups 
(\cite{Markman}).

On the dual side, 
let $f_1^+: X_1^+ \to X^+$ be the blowing up along the center $G^+= W^{+(0)}$, 
the $0$-section of the projection $\pi^+$.
Then the strict transform $(W^+)'$ of $W^+ = W^{+(1)}$ is smooth.
Let $f_2^+: X_2^+ \to X_1^+$ be the blowing up along the center $(W^+)'$.
Let $E^+_i$ be the exceptional divisor of $f_i^+$ for $i=1,2$, and 
$(E_1^+)'=(f_{2*}^+)^{-1}E_1^+$.
We set $f^+ = f_2^+ \circ f_1^+$.

The birational map $(\phi^+)^{-1} \circ \phi: X -\to X^+$ induces an 
isomorphism $X_2 \to X_2^+$ (\cite{Markman}).
We write $Y = X_2 = X_2^+$ under this identification.
Then we have $E_1' = (E_1^+)'$, $E_2 = E_2^+$, and
a commutative diagram
\[
\begin{CD}
A \in Z @<{\phi}<< X @<{f_1}<< X_1 @<{f_2}<< X_2 \\
@| @. @. @| \\
A^+ \in Z^+ @<{\phi^+}<< X^+ @<{f^+_1}<< X_1^+ @<{f^+_2}<< X_2^+
\end{CD}
\]
Since $\dim X = 4n-7$, $\dim G = 2n-4$ and $\dim W = 3n-5$, 
we have $K_Y = f^*K_X+(2n-4)E_1'+(n-3)E_2$.

The isomorphism $X_2 \to X_2^+$ can be 
described set theoretically in the following way.
 
(1) For a point $(p, A) \in X \setminus W$, we have 
$\text{rank }A = 2$ and $S_p = \text{Im }A$.
The corresponding point $(q, {}^tA) \in X^+ \setminus W^+$ is given by 
$S^+_q = \text{Im }{}^tA =(\text{Ker }A)^{\perp}$.

(2) For a point $\lim_{\epsilon \to 0} (p, \epsilon A_1) \in E_1 \setminus W'$,
with $\text{rank }A_1 = 2$ and $S_p = \text{Im }A_1$, the 
corresponding point is $\lim_{\epsilon \to 0} (q, \epsilon {}^tA_1) \in 
E_1^+ \setminus (W^+)'$, where $S^+_q = \text{Im }{}^tA_1$.

(3) For a point $(p, A) \in W \setminus G$, we have $\text{rank }A = 1$
and $S_p \supset \text{Im }A$. 
We take $A_2$ such that $S_p = \text{Im }(A+A_2)$.
Then a point 
$\lim_{\epsilon \to 0} (p, A + \epsilon A_2) \in E_2 \setminus E_1'$ 
corresponds to a point 
$\lim_{\epsilon \to 0} (q_{\epsilon}, {}^t(A + \epsilon A_2)) 
\in E_2^+ \setminus (E_1^+)'$ which is over a point $q \in G^+$ given by  
$S^+_q = \lim_{\epsilon \to 0} \text{Im }{}^t(A+\epsilon A_2)$.

(4) For a point $\lim_{\epsilon \to 0} (p, \epsilon A_1) \in E_1 \cap W'$ with 
$\text{rank }A_1 = 1$ and $S_p \supset \text{Im }A_1$, 
we take $A_2$ such that $S_p = \text{Im}(A_1+A_2)$.
Then a point 
$\lim_{\epsilon \to 0} (p, \epsilon A_1 + \epsilon^2 A_2) \in 
E_1' \cap E_2$ 
corresponds to a point 
$\lim_{\epsilon \to 0} (q_{\epsilon}, {}^t(\epsilon A_1 + \epsilon^2 A_2)) 
\in (E_1^+)' \cap E_2^+$ which is over a point $q \in G^+$ given by 
$S^+_q = \lim_{\epsilon \to 0} 
\text{Im }({}^t(A_1+\epsilon A_2))$.

The above description shows that the birational map $X_2 -\to X_2^+$
induced by $(\phi^+)^{-1} \circ \phi$ 
is an isomorphism by the Zariski main theorem.
As a consequence, we have the following lemma:

\begin{Lem}\label{SandS}
\[
S^*_X \vert_{X \setminus W} = S^+_{X^+} \vert_{X^+ \setminus W^+}.
\]
Moreover, $f^*S^*_X$ is a locally free subsheaf of $f^{+*}S^+_{X^+}$
such that a local section of the latter belongs to the former if and only if 
its value in $S^+_q$ at any point in $(f^+)^{-1}(x^+)$ 
for $x^+ = (q,B)$ is contained in $\text{Im }B$.
\end{Lem}

\begin{proof}
For $(p,A) \in X \setminus W$, we have 
$\text{Ker }{}^tA = (\text{Im }A)^{\perp}$ from (1) above.
Therefore, ${}^tA: V^* \to V^*$ induces an isomorphism 
$S_p^* = V^*/\text{Ker }{}^tA \to S^+_q$, 
hence the first assertion.

Any element $v \in V^*$ determines a global section of $f^*S^*_X$, and
the sheaf $f^*S^*_X$ is generated by such sections.
The value at the point $(q, {}^tA) \in X^+ \setminus W^+$ 
of the corresponding section of $S^+_{X^+}$ 
is given by ${}^tAv \subset S^+_q$, hence the second assertion.
\end{proof}

By taking the determinants, we obtain:

\begin{Cor}\label{OandO}
\[
f^*\mathcal{O}_X(1) 
= f^{+*}\mathcal{O}_{X^+}(-1) \otimes \mathcal{O}_Y(- 2E_1' - E_2).
\]
\hfill $\square$
\end{Cor}

Let $P$ be a fiber of $f_1$ above a point in $G$,
$P_0 = P \cap f_{1*}^{-1}X_0$ and $P_1 = P \cap W'$.  
Then the sequence $P_1 \subset P_0 \subset P$ is isomorphic to 
$\mathbb{P}^1 \times \mathbb{P}^{n-3} \subset \mathbb{P}^{2n-5} \subset 
\mathbb{P}^{2n-4}$, where the first inclusion is the Segre embedding and 
the second is linear.
Let $P' = f_{2*}^{-1}P$, $P'_0=f_{2*}^{-1}P_0$, and $P'_1=f_2^{-1}(P_1)$.
Let $l_1$ be a line on $P_0$ which meets $W'$ at $2$ points,
$l'_1 = f_{2*}^{-1}l_1$, and $l_2$ a fiber of $f_2$ over a point in $W'$.

\begin{Lem}
The cone of curves $\overline{NE}(Y/X)$ is generated by the classes of 
$l_1'$ and $l_2$.
The intersection numbers are given by the following table:
\[
(E_1' \cdot l'_1) = -1, 
(E_1' \cdot l_2) = 0, 
(E_2 \cdot l'_1) = 2, 
(E_2 \cdot l_2) = -1.
\]
\end{Lem}

\begin{proof}
$P_1$ is a determinantal variety defined by quadratic equations.
Therefore we obtain the formula for $\overline{NE}(Y/X)$.
The intersection numbers are obvious.
\end{proof}

%%%%%%%%%%%%%%%%%%%%%%%%%%%%%%%%%%%%%%%%%%%%%%%%%%%%%%%%%%%%%%%%%%%%

\section{Derived equivalence}

Let $X \to Z \leftarrow X^+$ and $X_0 \to Z_0 \leftarrow X^+_0$ 
be stratified flops defined in the previous section.
We shall compare the bounded derived categories of coherent sheaves
$D^b(\text{Coh}(X))$ with $D^b(\text{Coh}(X^+))$ and 
$D^b(\text{Coh}(X_0))$ with $D^b(\text{Coh}(X^+_0))$.
By \cite{Kawamata2}~Lemma~5.6 and Corollary~5.7, the derived equivalence for 
the former pair implies the latter.

We begin with an easy case $r = 1$, i.e., the standard flop and the 
usual Mukai flop.

\begin{Prop}
Let $l$ be an arbitrary integer.
Then a functor 
\[
\Phi_l: D^b(\text{Coh}(X)) \to D^b(\text{Coh}(X^+))
\]
defined by
\[
\Phi_l(a) = f^+_*(f^*a \otimes \mathcal{O}_Y(lE))
\]
is an equivalence.
\end{Prop}
 
\begin{proof}
The adjoint functor $\Psi_l: D^b(\text{Coh}(X^+)) \to D^b(\text{Coh}(X))$
is given by 
\[
\Psi_l(b) = f_*(f^{+*}a \otimes \mathcal{O}_Y((n-1-l)E)).
\]
The category $D^b(\text{Coh}(X))$ is spanned by the set of sheaves
$\mathcal{O}_X(k)$ for $l - n + 1 \le k \le l$.
By the Kodaira vanishing theorem, we have
\[
\Phi_l: \mathcal{O}_X(k) \mapsto \mathcal{O}_Y(k,0)(lE)
= \mathcal{O}_Y(0,-k)((l-k)E)
\mapsto \mathcal{O}_{X^+}(-k)
\]
and
\[
\begin{split}
\Psi_l: &\mathcal{O}_{X^+}(-k) \mapsto \mathcal{O}_Y(0,-k)((n-1-l)E)
= \mathcal{O}_Y(k,0)((n-1-l+k)E) \\
&\mapsto \mathcal{O}_X(k).
\end{split}
\]
There is an adjunction morphism of functors
$F: \text{Id}_{D^b(\text{Coh}(X))} \to \Psi_l\Phi_l$,
which is reduced to the identity when restricted to the open subset 
$X \setminus W$.
By the above argument, we have isomorphisms 
$\omega \cong \Psi\Phi(\omega)$ for a spanning class $\Omega = \{\omega\}$.
Since the $\omega$ are invertible sheaves, 
it follows that the morphisms $F(\omega)$ are isomorphisms.
Therefore, the natural homomorphisms
\[
\Phi: \text{Hom}^p(\omega_1, \omega^2) 
\to \text{Hom}^p(\omega_1, \Psi\Phi(\omega^2)) 
\to \text{Hom}^p(\Phi(\omega_1), \Phi(\omega^2))  
\]
for any $\omega_1, \omega_2 \in \Omega$ and $p \in \mathbb{Z}$ 
are isomorphisms.
Then $\Phi$ is an equivalence by \cite{BO} and \cite{Bridgeland1}.
\end{proof}

We assume that $r = 2$ in the rest of the paper.
We consider exact functors between bounded derived categories
\begin{equation}\label{PhiPsi}
\begin{split}
&\Phi: D^b(\text{Coh}(X)) \to D^b(\text{Coh}(Y)) \\
&\Psi: D^b(\text{Coh}(Y)) \to D^b(\text{Coh}(X))
\end{split}
\end{equation}
defined by
\[
\begin{split}
&\Phi(a) = Rf^+_*(Lf^*(a) \otimes \mathcal{O}_Y((2n-5)E_1'+(n-3)E_2)) \\
&\Psi(b) = Rf_*(Lf^{+*}(b) \otimes \mathcal{O}_Y(E_1')).
\end{split}
\]
They are adjoints each other because $K_Y = f^*K_X+(2n-4)E_1'+(n-3)E_2$.

\begin{Lem}\label{span}
$D^b(\text{Coh}(X))$ is spanned by the set of the following 
locally free sheaves:
\[
\text{Sym}^iS^*_X \otimes \mathcal{O}_X(j)
\]
for $0 \le i$, $0 \le j$ and $i+j\le n-2$.
\end{Lem}

\begin{proof}
By \cite{Kapranov}, any point sheaf $\mathcal{O}_p$ for $p \in G$ has a finite
locally free resolution whose terms are direct sums of the sheaves 
\[
\text{Sym}^iS^* \otimes \mathcal{O}_G(j)
\] 
for $0 \le i$, $0 \le j$ and $i+j\le n-2$.
Hence $\mathcal{O}_{\pi^{-1}(p)}$ is resolved by our set.
Then so is any point sheaf $\mathcal{O}_x$ for $x \in X$
because $\pi^{-1}(p)$ is an affine space.
\end{proof}

Let $\mathcal{E}^+$ be a subsheaf of $S^+_{X^+}$
such that a local section of the latter belongs to the former if and only if 
its value in $S^+_q$ at a point $(q,B)$ is contained in $\text{Im }B$.
We denote by $\mathcal{E}^+_i$ for $i > 0$ 
the image of the natural homomorphism
$\text{Sym}^{i-1}S^+_{X^+} \otimes \mathcal{E}^+ \to \text{Sym}^iS^+_{X^+}$.

\begin{Lem}\label{PhiPsi2}
(1) For $0 \le j \le n-3$,
\[
\begin{split}
&\Phi(\mathcal{O}_X(j))=\mathcal{O}_{X^+}(-j) \\
&\Psi(\mathcal{O}_{X^+}(-j))=\mathcal{O}_X(j).
\end{split}
\]

(2) Let $I_{W^+}$ be the ideal sheaf of $W^+ \subset X^+$.  Then 
\[
\Phi(\mathcal{O}_X(n-2))=
I_{W^+} \otimes \mathcal{O}_{X^+}(-n+2).
\]

(3) For $0 <i$, $0 \le j$ and $i+j \le n-3$,
\[
\begin{split}
&\Phi(\text{Sym}^iS^*_X \otimes \mathcal{O}_X(j))=
\text{Sym}^iS^+_{X^+} \otimes \mathcal{O}_{X^+}(-j) \\
&\Psi(\text{Sym}^iS^+_{X^+} \otimes \mathcal{O}_{X^+}(-j))=
\text{Sym}^iS^*_X \otimes \mathcal{O}_X(j).
\end{split}
\]

(4) For $0 < i \le n-2$,
\[
\Phi(\text{Sym}^iS^*_X \otimes \mathcal{O}_X(n-2-i))=
\mathcal{E}^+_i \otimes \mathcal{O}_{X^+}(-n+2+i).
\]
\end{Lem}

\begin{proof}
(1) By Corollary~\ref{OandO} 
\[
\begin{split}
&f^*(\mathcal{O}_X(j)) \otimes \mathcal{O}_Y((2n-5)E_1'+(n-3)E_2) \\
&= f^{+*}\mathcal{O}_{X^+}(-j) \otimes 
\mathcal{O}_Y((2n-5-2j)E_1'+(n-3-j)E_2).
\end{split}
\]
Since $K_Y = f^*K_X+(2n-4)E_1'+(n-3)E_2$,
the higher direct images for $f^+$ vanish because 
$-(2j+1)E_1'-jE_2$ is nef.
Since $(2n-5-2j)E_1'+(n-3-j)E_2$ is effective, we have
the first formula. 

Similarly, we have
\[
f^{+*}(\mathcal{O}_{X^+}(-j)) \otimes \mathcal{O}_Y(E_1')
= f^*\mathcal{O}_X(j) \otimes \mathcal{O}_Y((2j+1)E_1'+jE_2).
\]
The higher direct images for $f$ vanish because 
$(2j-2n+5)E_1'+(j-n+3)E_2$ is nef.
$(2j+1)E_1'+jE_2$ is effective, hence
the second formula.

(2) We already proved in (1) that the higher direct images vanish.
We obtain our formula from 
\[
f^+_*\mathcal{O}_Y(-E_1'-E_2) = I_{W^+}.
\]

(3) We have an exact sequence
\[
0 \to \mathcal{O}_{E_2}(- E_1') \to 
f^{+*}S^+_{X^+} \otimes \mathcal{O}_{E_2}(- E_1')
\to \mathcal{O}_{E_2}(- E_1') \otimes f^{+*}\mathcal{O}_{X^+}(-1) \to 0.
\]
where the first term is the image of $f^*S^*_X$.
We define a decreasing filtration of the sheaf
$f^*\text{Sym}^iS^*_X$
by locally free subsheaves $\mathcal{F}^{k,l}$ for $0 \le l \le k \le i$ by
\[
\begin{split}
\mathcal{F}^{k,l} 
= &f^*\text{Sym}^iS^*_X \cap f^{+*}\text{Sym}^iS^+_{X^+}(-iE_1'-kE_2) \\
&\cap (f^*\text{Sym}^iS^*_X(-lE_2)
+f^{+*}\text{Sym}^iS^+_{X^+}(-iE_1'-(k+1)E_2)). 
\end{split}
\]
We have
\[
\begin{split}
&f^*\text{Sym}^iS^*_X = F^{0,0} \supset F^{1,0} \supset F^{1,1} \supset 
F^{2,0} \supset \\
&\dots \supset F^{i-1,i-1} \supset 
F^{i,0}=f^{+*}\text{Sym}^iS^+_{X^+}(-iE_1'-iE_2).
\end{split}
\]
and 
\[
\mathcal{F}^{k,l}/\mathcal{F}^{k,l+1} 
\cong \mathcal{O}_{E_2}(- iE_1'-kE_2) 
\otimes f^{+*}\mathcal{O}_{X^+}(-k+l)
\]
for $0 \le l \le k < i$, 
where we put $\mathcal{F}^{k,k+1}=\mathcal{F}^{k+1,0}$.
Since 
\[
R(f_2^+)_*\mathcal{O}_{E_2}(tE_2)=0
\]
for $0 < t \le n-3$, 
we have
\[
R(f_2^+)_*(\mathcal{F}^{k,l}/\mathcal{F}^{k,l+1}
\otimes \mathcal{O}_Y((2n-5-2j)E'_1+(n-3-j)E_2)) = 0 
\]
because $i + j \le n-3$.
Therefore,
\[
\begin{split}
\Phi(&\text{Sym}^iS^*_X \otimes \mathcal{O}_X(j)) \\
=&\text{Sym}^iS^+_{X^+} \otimes \mathcal{O}_{X^+}(-j) \\
&\otimes Rf^+_*\mathcal{O}_Y((2n-5-i-2j)E'_1+(n-3-i-j)E_2) \\
=&\text{Sym}^iS^+_{X^+} \otimes \mathcal{O}_{X^+}(-j)
\otimes Rf^+_{1*}\mathcal{O}_{X_1^+}((2n-5-i-2j)E_1^+) \\
=&\text{Sym}^iS^+_{X^+} \otimes \mathcal{O}_{X^+}(-j).
\end{split}
\]

For the inverse direction, we climb back along the filtration 
$\mathcal{F}^{k,l}$.
Since 
\[
Rf_{2*}(\mathcal{F}^{k,l}/\mathcal{F}^{k,l+1}
\otimes \mathcal{O}_Y((i+2j+1)E'_1+(i+j)E_2)) = 0 
\]
for $0 \le l \le k < i$, we have
\[
\begin{split}
&\Psi(\text{Sym}^iS^+_{X^+} \otimes \mathcal{O}_{X^+}(-j)) \\
&=\text{Sym}^iS^*_X \otimes \mathcal{O}_X(j)
\otimes Rf_*\mathcal{O}_Y((i+2j+1)E'_1+(i+j)E_2) \\
&=\text{Sym}^iS^*_X \otimes \mathcal{O}_X(j)
\otimes Rf_{1*}\mathcal{O}_{X_1}((i+2j+1)E_1) \\
&=\text{Sym}^iS^*_X \otimes \mathcal{O}_X(j).
\end{split}
\]

(4) We already proved in (3) that
\[
\begin{split}
&Rf^+_*(f^*(\text{Sym}^iS^*_X \otimes \mathcal{O}_X(n-2-i)) 
\otimes \mathcal{O}_Y((2n-5)E_1'+(n-3)E_2)) \\
&\cong Rf^+_*(F^{i-1,0} \otimes f^{+*}\mathcal{O}_X(-n+2+i)
\otimes \mathcal{O}_Y((2i-1)E_1'+(i-1)E_2)).
\end{split}
\]
Moreover
\[
\begin{split}
&R^pf^+_*(f^{+*}(\text{Sym}^iS^+_{X^+} \otimes \mathcal{O}_{X^+}(-n+2+i))
\otimes \mathcal{O}_Y((i-1)E_1'-E_2)) \\
&\to R^pf^+_*(F^{i-1,0} \otimes f^{+*}\mathcal{O}_X(-n+2+i)
\otimes \mathcal{O}_Y((2i-1)E_1'+(i-1)E_2)) 
\end{split}
\]
is surjective for $p > 0$.
Therefore, there are no higher direct images for $f^+$.
Since
\[
\begin{split}
&0 \to F^{i-1,0} 
\to f^{+*}\text{Sym}^iS^+_{X^+} \otimes \mathcal{O}_Y(-iE_1'-(i-1)E_2) \\
&\to \mathcal{O}_{E_2}(-iE_1'-(i-1)E_2) \otimes f^{+*}\mathcal{O}_{X^+}(-i)
\to 0
\end{split}
\]
we have 
\[
\begin{split}
&0 \to f^+_*(F^{i-1,0} \otimes f^{+*}\mathcal{O}_X(-n+2+i)
\otimes \mathcal{O}_Y((2i-1)E_1'+(i-1)E_2)) \\
&\to \text{Sym}^iS^+_{X^+} \otimes \mathcal{O}_X(-n+2+i)
\to \mathcal{O}_{W^+} \otimes \mathcal{O}_{X^+}(-n+2)
\to 0.
\end{split}
\]
On the other hand, we have 
\[
\begin{split}
&f^+_*(F^{i-1,0} \otimes f^{+*}\mathcal{O}_X(-n+2+i)
\otimes \mathcal{O}_Y((2i-1)E_1'+(i-1)E_2)) \vert_{X^+ \setminus G^+} \\
&\cong \mathcal{E}^+_i \otimes \mathcal{O}_{X^+}(-n+2+i) 
\vert_{X^+ \setminus G^+}.
\end{split}
\]
By taking the direct image sheaves of both sides under the inclusion morphism
$X^+ \setminus G^+ \to X^+$, we obtain 
\[
\begin{split}
&f^+_*(F^{i-1,0} \otimes f^{+*}\mathcal{O}_X(-n+2+i)
\otimes \mathcal{O}_Y((2i-1)E_1'+(i-1)E_2)) \\
&\cong \mathcal{E}^+_i \otimes \mathcal{O}_{X^+}(-n+2+i).
\end{split}
\]
\end{proof}

We assume that $n = 4$ besides $r = 2$ in the following.
Then $\dim G = 4$, $\dim X_0 = 8$, 
$\dim X = 9$, and $\dim W = 7$.
In particular, $W$ is locally complete intersection.
We have $K_Y = f^*K_X+4E_1'+E_2$.

In this case, Namikawa \cite{Namikawa2} proved that the functors
\[
\Phi': D^b(\text{Coh}(X)) \to D^b(\text{Coh}(Y)), 
\Psi': D^b(\text{Coh}(Y)) \to D^b(\text{Coh}(X))
\]
defined by
\[
\Phi'(a) = Rf^+_*(Lf^*(a) \otimes \mathcal{O}_Y(4E_1'+E_2)),
\Psi'(b) = Rf_*(Lf^{+*}(b) \otimes \mathcal{O}_Y).
\]
are {\em not} equivalences.

\begin{Lem}\label{PhiPsi3}
Assume that $n=4$. Then

(1) $\Psi(I_{W^+} \otimes \mathcal{O}_{X^+}(-2))
=\mathcal{O}_X(2)$.

(2) $\Psi(\mathcal{E}^+_i \otimes \mathcal{O}_{X^+}(i-2))=
\text{Sym}^iS^* \otimes \mathcal{O}_X(2-i)$ 
for $0 < i \le 2$.
\end{Lem}

\begin{proof}
(1) $W^+$ is a divisor on $X_0^+$ and $\phi^+$ is the contraction 
of $(-2)$-curves along the generic points of $\phi(W^+)$.
Hence $\mathcal{O}_{X_0^+}(-W^+) = \mathcal{O}_{X_0^+}(2)$.
We have an exact sequence
\[
0 \to \mathcal{O}_{X^+}(-2-X_0^+) \to I_{W^+} \otimes \mathcal{O}_{X^+}(-2) 
\to \mathcal{O}_{X_0^+} \to 0.
\]
Since $X_0^+$ is a Cartier divisor, we have 
$L_kf^{+*}\mathcal{O}_{X_0^+} = 0$ for $k > 0$. 
Hence
\[
0 \to f^{+*}\mathcal{O}_{X^+}(-2-X_0^+) \to 
f^{+*}(I_{W^+} \otimes \mathcal{O}_{X^+}(-2)) \to 
f^{+*}\mathcal{O}_{X_0^+} \to 0.
\]
On the other hand,
\[
0 \to f^{+*}\mathcal{O}_{X^+}(-2-X_0^+) \to 
I_{W^+}\mathcal{O}_Y \otimes f^{+*}\mathcal{O}_{X^+}(-2) \to 
\mathcal{O}_{(X_0^+)''} \to 0
\]
where $(X_0^+)'' = (f^+_*)^{-1}X_0^+$.
Thus
\[
\begin{split}
&0 \to \mathcal{O}_{(E_1^+)' + E_2^+}(- (X_0^+)'') \to 
f^{+*}(I_{W^+} \otimes \mathcal{O}_{X^+}(-2)) \\
&\to I_{W^+}\mathcal{O}_Y \otimes f^{+*}\mathcal{O}_{X^+}(-2) \to 0.
\end{split}
\]
We have
\[
\mathcal{O}_{(E_1^+)' + E_2^+}(-(X_0^+)'')
\cong \mathcal{O}_{E_1' + E_2}(E_1'+E_2).
\]
Since $Rf_{2*}\mathcal{O}_{E_2}(E_2)=0$ and $Rf_{1*}\mathcal{O}_{E_1}(2E_1)=0$,
we have
\[
Rf_*\mathcal{O}_{E_1' + E_2}(2E_1'+E_2) 
\cong Rf_{1*}\mathcal{O}_{E_1}(2E_1) \cong 0.
\] 
Therefore, we have
\[
\begin{split}
&Rf_*(f^{+*}(I_{W^+} \otimes \mathcal{O}_{X^+}(-2)) 
\otimes \mathcal{O}_Y(E_1')) \\
&\cong Rf_*(I_{W^+}\mathcal{O}_Y \otimes f^{+*}\mathcal{O}_{X^+}(-2)
\otimes \mathcal{O}_Y(E_1')).
\end{split}
\]

Let $F^+ = (X^+_0)'' \cap (E^+_1)' \subset Y$. 
The projection $f^+: F^+ \to G^+$ is a $\mathbb{P}^3$-bundle whose fiber 
$(P_0^+)'$ is contracted by $f$.
We have 
\[
I_{W^+}\mathcal{O}_Y = I_{F^+}(- E_1' - E_2)
\]
Thus
\[
I_{W^+}\mathcal{O}_Y \otimes f^{+*}\mathcal{O}_{X^+}(-2) 
\otimes \mathcal{O}_Y(E_1')
= I_{F^+}f^*\mathcal{O}_X(2) \otimes \mathcal{O}_Y(4E_1'+E_2).
\]
On the other hand, there is an exact sequence
we have $\mathcal{O}_{(P_0^+)'}(E_1') \cong \mathcal{O}_{\mathbb{P}^3}(-1)$ 
and $\mathcal{O}_{(P_0^+)'}(E_2) \cong \mathcal{O}_{\mathbb{P}^3}(2)$.
Hence $\mathcal{O}_{(P_0^+)'}(4E_1'+E_2) \cong \mathcal{O}_{\mathbb{P}^3}(-2)$.
It follows that $Rf_*\mathcal{O}_{F^+}(4E_1'+E_2)=0$.
Therefore, 
\[
Rf_*(I_{F^+}\mathcal{O}_Y \otimes f^{+*}\mathcal{O}_{X^+}(-2) 
\otimes \mathcal{O}_Y(E_1')) \cong \mathcal{O}_X(2).
\]

(2) We have an exact sequence
\[
0 \to \mathcal{E}^+_i \otimes \mathcal{O}_{X^+}(i-2) 
\to \text{Sym}^iS^+_{X^+} \otimes \mathcal{O}_{X^+}(i-2)
\to \mathcal{O}_{W^+} \otimes \mathcal{O}_{X^+}(-2) \to 0.
\]
Thus the torsion part of 
$f^{+*}\mathcal{E}^+_i \otimes \mathcal{O}_{X^+}(i-2)$ is isomorphic to that 
of $f^{+*}I_{W^+} \otimes \mathcal{O}_{X^+}(-2)$.
Therefore, by the argument in (1), we have
\[
\Psi(\mathcal{E}^+_i \otimes \mathcal{O}_{X^+}(i-2))
\cong Rf_*(f^{+*}\mathcal{E}^+_i/\text{torsion} \otimes 
f^{+*}\mathcal{O}_{X^+}(i-2) \otimes \mathcal{O}_Y(E'_1)).
\]
We have an exact sequence
\[
\begin{split}
&0 \to f^{+*}\mathcal{E}^+_i/\text{torsion} 
\otimes f^{+*}\mathcal{O}_{X^+}(i-2) \otimes \mathcal{O}_Y(E'_1) \\
&\to f^{+*}(\text{Sym}^iS^+_{X^+} \otimes \mathcal{O}_{X^+}(i-2))
\otimes \mathcal{O}_Y(E'_1) \\
&\to f^{+*}(\mathcal{O}_{W^+} \otimes \mathcal{O}_{X^+}(-2)) 
\otimes \mathcal{O}_Y(E'_1) \to 0.
\end{split}
\]
Therefore, we have an exact sequence
\[
\begin{split}
&0 \to f^{+*}\mathcal{E}^+_i/\text{torsion} 
\otimes f^{+*}\mathcal{O}_{X^+}(i-2) \otimes \mathcal{O}_Y(E'_1) \\
&\to F^{i-1,0} \otimes f^{+*} \mathcal{O}_{X^+}(i-2)
\otimes \mathcal{O}_Y((i+1)E'_1+(i-1)E_2)
\to Q \to 0
\end{split}
\]
where the cokernel $Q$ has a decomposition as follows
\[
0 \to \mathcal{O}_{F^+}(- E_2) \otimes f^{+*}\mathcal{O}_{X^+}(-2)
\to Q \to \mathcal{O}_{E'_1}(E'_1 - E_2) \otimes f^{+*}\mathcal{O}_{X^+}(-2)
\to 0.
\]
As in (1), $F^+=F=X_0'' \cap E_1'$ is a $\mathbb{P}^3$-bundle over $G$ with
fibers $P_0'$ such that $\mathcal{O}_{P_0'}(E_1') \cong 
\mathcal{O}_{\mathbb{P}^3}(-1)$ and $\mathcal{O}_{P_0'}(E_2)=
\mathcal{O}_{\mathbb{P}^3}(2)$.
Thus
\[
\mathcal{O}_{F^+}(- E_2) \otimes f^{+*}\mathcal{O}_{X^+}(-2) \otimes
\mathcal{O}_{P_0'} \cong \mathcal{O}_{P_0'}(4E'_1+E_2) 
\cong \mathcal{O}_{P_0'}(-2)
\]
and we have $Rf_*\mathcal{O}_{F^+}(- E_2)=0$.
We have also
\[
\begin{split}
&Rf_*\mathcal{O}_{E'_1}(E'_1 - E_2) \otimes f^{+*}\mathcal{O}_{X^+}(-2) \\
&\cong Rf_*\mathcal{O}_{E'_1}(3E'_1+E_2) \otimes \mathcal{O}_X(2) \\
&\cong Rf_{1*}\mathcal{O}_{E_1}(3E_1) \otimes \mathcal{O}_X(2)
\cong 0.
\end{split}
\]
Therefore
\[
\begin{split}
&\Psi(\mathcal{E}^+_i \otimes \mathcal{O}_{X^+}(i-2)) \\
&\cong Rf_*(F^{i-1,0} \otimes f^{+*} \mathcal{O}_{X^+}(i-2)
\otimes \mathcal{O}_Y((i+1)E'_1+(i-1)E_2)) \\
&\cong Rf_*(F^{i-1,0} \otimes \mathcal{O}_Y((5-i)E'_1+E_2)) 
\otimes \mathcal{O}_X(2-i).
\end{split}
\]
If $i=1$, then $F^{0,0}=f^*\text{Sym}^iS^*$ and
$Rf_*\mathcal{O}_Y(4E'_1+E_2) = \mathcal{O}_X$, hence the result.

If $i=2$, then $\mathcal{F}^{0,0}/\mathcal{F}^{1,0} \cong 
\mathcal{O}_{E_2}(-2E_1')$, hence
\[
Rf_*(\mathcal{F}^{0,0}/\mathcal{F}^{1,0} \otimes \mathcal{O}_Y(3E'_1+E_2))
= 0
\]
and we complete the proof.
\end{proof}

\begin{Thm}
If $n=4$, then $\Phi$ and $\Psi$ in (\ref{PhiPsi}) are equivalences.
\end{Thm}

\begin{proof}
There is an adjunction morphism of functors
$F: \text{Id}_{D^b(\text{Coh}(X))} \to \Psi\Phi$,
which is reduced to the identity when restricted to the open subset 
$X \setminus W$.
By Lemmas~\ref{span}, \ref{PhiPsi2} and \ref{PhiPsi3}, we have isomorphisms
$\omega \cong \Psi\Phi(\omega)$ for the spanning class $\Omega = \{\omega\}$ 
given by Lemma~\ref{span}.
Since the $\omega$ are locally free sheaves, 
it follows that the morphisms $F(\omega)$ are isomorphisms.
Therefore, the natural homomorphisms
\[
\Phi: \text{Hom}^p(\omega_1, \omega^2) 
\to \text{Hom}^p(\omega_1, \Psi\Phi(\omega^2)) 
\to \text{Hom}^p(\Phi(\omega_1), \Phi(\omega^2))  
\]
for any $\omega_1, \omega_2 \in \Omega$ and $p \in \mathbb{Z}$ 
are isomorphisms.
Then $\Phi$ is an equivalence by \cite{BO} and \cite{Bridgeland1}.
\end{proof}

\begin{Rem}
(1) Let $F = \bigoplus \omega$ be the sum of the spanning sheaves 
given in Lemma~\ref{span}.
Then $F$ is {\em not} an almost exceptional object in the sense of \cite{BK}
in our case $n=4$ and $r=2$,
because $R^1\phi_*\mathcal{O}_X(-2) \ne 0$ implies $\text{Hom}^1(F,F) \ne 0$.

(2) In order to extend our argument to the case $n \ge 5$, 
the {\em Eagon-Northcott resolution} on $X_0^+$ 
would be useful (\cite{Eisenbud}):
\[
\begin{split}
&0 \to \pi_0^{+*}(\text{Sym}^{n-2}S^{+2})^* \otimes \bigwedge^n V^{*n} 
\stackrel{B_*}{\to} \dots 
\to \pi_0^{+*}(\text{Sym}^kS^{+2})^* \otimes \bigwedge^{k+2} V^{*n} \\
&\stackrel{B_*}{\to} 
\pi_0^{+*}(\text{Sym}^{k-1}S^{+2})^* \otimes \bigwedge^{k+1} V^{*n} 
\stackrel{B_*}{\to} \dots 
\stackrel{B_*}{\to} \pi_0^{+*}(S^{+2})^* \otimes \bigwedge^3 V^{*n} \\
&\stackrel{B_*}{\to} \mathcal{O}_{X^+_0} \otimes \bigwedge^2 V^{*n} 
\stackrel{\bigwedge^2 B}{\to} \pi_0^{+*} \bigwedge^2 S^{+2} \to 
\mathcal{O}_{W^+}(-1) \to 0.
\end{split}
\]
\end{Rem}

%%%%%%%%%%%%%%%%%%%%%%%%%%%%%%%%%%%%%%%%%%%%%%%%%%%%%%%%%%%%%%%%%%%%

Department of Mathematical Sciences, University of Tokyo, 

Komaba, Meguro, Tokyo, 153-8914, Japan 

kawamata@ms.u-tokyo.ac.jp

\end{document}